\documentclass[11pt]{article}

\input{amssym.def}
\input{amssym}
\usepackage{eufrak}
\usepackage{color}

\usepackage{graphicx}

\newtheorem{theorem}{Theorem}

\newtheorem{corollary}[theorem]{Corollary}
\newtheorem{definition}[theorem]{Definition}

\newtheorem{proposition}[theorem]{Proposition}
\newtheorem{remark}[theorem]{Remark}

\def\pa{{\partial}}

\def\ov{{\overline}}
\def\Y{\mathbf{Y}}

\def\qq{q^{-1}}

\def\Tr{\rm Tr}

\def\A{{\cal A}}

\def\FF{{\cal F}}

\def\XX{\frak{X}}

\def\ot{\otimes}

\def\C{{\Bbb C}}

\def\vv{V^{\otimes 2}}
\def\RR{{\cal R}}

\def\ov{\overline}

\def\be{\begin{equation}}
\def\ee{\end{equation}}

\def\bea{\begin{eqnarray}}
\def\eea{\end{eqnarray}}
\def\nn{\nonumber}

\begin{document}

\makeatletter
\renewcommand{\theequation}{{\thesection}.{\arabic{equation}}}
\@addtoreset{equation}{section} \makeatother

\title{KZ equations and Bethe subalgebras in generalized Yangians related to compatible $R$-matrices}

\author{\rule{0pt}{7mm} Dimitri
Gurevich\thanks{gurevich@ihes.fr}\\
{\small\it Univ. de Valenciennes, EA 4015-LAMAV}\\
{\small\it F-59313 Valenciennes, France}\\
{\small \it and}\\
{\small \it Interdisciplinary Scientific Center J.-V.Poncelet}\\
{\small\it Moscow 119002, Russian Federation}\\
\rule{0pt}{7mm} Pavel Saponov\thanks{Pavel.Saponov@ihep.ru}\\
{\small\it
National Research University Higher School of Economics,}\\
{\small\it 20 Myasnitskaya Ulitsa, Moscow 101000, Russian Federation}\\
{\small \it and}\\
{\small \it
Institute for High Energy Physics, NRC "Kurchatov Institute"}\\
{\small \it Protvino 142281, Russian Federation}\\
\rule{0pt}{7mm} Dmitry Talalaev\thanks{dtalalaev@yandex.ru}\\
{\small\it
Moscow State University, Faculty of Mechanics and Mathematics}\\
{\small\it 119991 Moscow, Russian Federation}
\\
{\small \it and}\\
{\small \it
Institute for Theoretical and Experimental Physics, NRC "Kurchatov Institute"}\\
{\small \it 117218 Moscow, Russian Federation, 25 Bolshaya Cheremushkinskaya str.}\\
{\small \it and}\\
{\small \it
Centre of integrable systems, P. G. Demidov Yaroslavl State University}\\
{\small \it 150003, Yaroslavl, Russian Federation, 14 Sovetskaya str.}}

\maketitle

\begin{abstract}
The notion of compatible braidings was introduced in \cite{IOP}. On the base of  this notion the authors of \cite{IOP} defined certain quantum matrix algebras  generalizing the RTT algebras and Reflection Equation ones. They also defined  analogs of some symmetric polynomials in these algebras and showed that these polynomials  generate  commutative subalgebras, called Bethe. By using a similar approach we introduce certain new algebras called
generalized Yangians and  define analogs of some symmetric polynomials in these algebras. We claim that they commute with each other and thus generate a commutative
Bethe subalgebra in each generalized Yangian.
Besides, we define some analogs (also arising from couples of compatible braidings) of the Knizhnik-Zamolodchikov  equation--classical and quantum.
\end{abstract}

{\bf AMS Mathematics Subject Classification, 2010:} 81R50

{\bf Keywords:} compatible braidings, braided Yangians, Bethe subalgebra, braided $r$-matrix, braided (quantum) KZ connections

\section{Introduction}

The notion of compatible $R$-matrices (we call them {\em braidings})  was introduced in \cite{IOP}. By a braiding we mean a linear operator $R:\vv\to\vv$ subject to
the braid relation
\be (R\ot I)(I\ot R)(R\ot I)=(I\ot R)(R\ot I)(I\ot R), \label{braid} \ee
where  $V,\,\, dim V=N<\infty$, is a vector space over the ground field $\C$ and $I$ is the identity operator or its matrix.

According to \cite{IOP}, two braidings  $R$ and $F$ are called {\em compatible}   if they are subject to the  system
\be  R_{12}\, F_{23}\, F_{12} = F_{23}\, F_{12}\, R_{23},\,\,\,R_{23}\, F_{12}\, F_{23} = F_{12}\, F_{23}\, R_{12}.  \label{coor} \ee

As usual, the low indexes   indicate the positions, where a matrix (or an operator) is located.
Observe that  the matrices $A_k, \, k\geq 2$ are obtained from that $A_1$ by means of the usual flip  $P$ :
$$A_2=P_{12}\, A_1\, P_{12},   \, A_3=P_{23}\, A_2\, P_{23}=P_{23}\, P_{12}\, A_1\, P_{12}\, P_{23},\,\,{\rm  and\,\, so\,\, on}.$$
Here,  $A$ is an $N\times N$ matrix, whereas $A_1=A\ot I_{2...p},\,\, A_2=I_1\ot A\ot I_{3...p}$ and so on  are  $N^p\times N^p$ matrices.

Following  \cite{IOP}  we introduce the following notations
$$A_{\ov 2}=F_{12}\, A_1\, F_{12}^{-1},\,\, A_{\ov 3}=F_{23}\, A_{\ov 2}\, F_{23}^{-1}=F_{23}\,F_{12}\, A_1\,F_{12}^{-1}\, F_{23}^{-1},$$
and so on, where the   overlined indexes mean that the matrix $A_1$ is pushed forward to higher positions by means of   the second braiding $F$.
For the sake of the uniformity we also put $A_{\ov 1}=A_1$.

Below, we fix a basis in the space $V$ and the corresponding bases in $V^{\ot p}$  and  identify operators and their matrices.

Following \cite{IOP} introduce  an algebra $\A(R,F)$ defined by the following system of  relations
\be R_{12}\, T_{\ov 1}\,T_{\ov 2}=T_{\ov 1}\,T_{\ov 2}\,R_{12}, \label{RTT}\ee
where $T=\|t_i^j\|_{1\leq i,j\leq N}$ is a matrix with  entries $t_i^j$.

The matrix $T$ is called {\em generating matrix} of the algebra $\A(R,F)$.

As shown in \cite{IOP}, in these algebras (under some conditions on the braidings) it is possible to define analogs of some symmetric polynomials and to establish analogs of the Cayley-Hamilton and Newton identities. (In \cite{IOP} these identities are combined  in the so-called Cayley-Hamilton-Newton ones.) Also, it is possible to
show that {\em quantum elementary symmetric polynomials} commute with each other and consequently generate a subalgebra  of $\A(R,F)$ called {\em Bethe}.

The first purpose of the present paper is to introduce some algebras similar to those $\A(R,F)$ but with infinite number of generators and to generalize the mentioned results to them.
Each of these algebras is defined via the system
\be  R(u,v)\, T_{\ov{1}}(u)\, T_{\ov{2}}(v)=T_{\ov{1}}(v)\, T_{\ov{2}}(u)\, R(u,v), \label{brY} \ee
where, $T(u)=\sum_{k\geq 0} T[k]u^{-k}$ is a matrix  expanded in a Laurent series
 and the  current (i.e. depending on parameters) quantum  $R$-matrices $R(u,v)$ arises  from
a braiding $R$ via the Baxterization procedure. We denote this algebra $\Y(R,F)$ and call it the {\em generalized  Yangian}\footnote{Note that if
$$R(u,v)=P-\frac{I}{u-v}\,\,{\rm  and}\,\, F=P,$$ we get the famous Drinfeld's Yangian $Y(gl(N))$. (Usually, one also imposes the condition $T[0]=I$.)
Its generators are entries of the matrices 
$$T[k]=\|t_i^j[k]\|_{1\leq i,j \leq N},\,\, k=0,1,2...$$
A similar treatment  is valid for our generalized Yangians.
However, below we do not use this treatment  and deal  with the {\em generating}  matrix
$T(u)$ in whole.}.

So far, our generalized Yangians  are the most general quantum matrix algebras associated with rational and trigonometric $R$-matrices for which analogs of some symmetric polynomials,  namely, elementary ones and power sums, are constructed and their commutativity is established.

The second purpose of the paper is to introduce braided analogs of the Knizhnik-Zamolodchi\-kov (KZ) equation--classical and quantum--and
to establish their compatibility.    Observe that the former ones are based on a braided version of the first Sklyanin bracket. In its turn this version is based on  braided
 current $r$-matrices. In the rational case such a braided   current $r$-matrix can be easily defined with the help of an involutive symmetry $F$ as follows
$r(u,v)=\frac{F}{u-v}$. (So, in this case any second braiding is not needed.) 

In the trigonometric case we are looking for a braided current $r$-matrix under the form $r(u,v)=\frac{F\, u}{u-v}+r$. In order to find a (constant) summand $r$ we need a Hecke symmetry $R=R(q)$  analytically depending  on $q$ in a vicinity of $q=1$  and deforming
an involutive symmetry $F$ (i.e. $R(1)=F$).  Then by expanding $\RR=\RR(q)=R(q)\, F$ at the point $q=1$, we get  $r$.

The paper is organized as follows. In the next section we define compatible braidings and exhibit some examples. In Section 3 we introduce the aforementioned quantum
symmetric polynomials in the generalized Yangians. In Section 4 we describe braided versions of current $r$-matrices and the first Sklyanin bracket.
In section 5 we introduce braided analogs of the classical and quantum KZ equation in the spirit of \cite{KZ} and \cite{FR} respectively.

 {\bf Acknowledgements:} The work of P.S. has been funded by the Russian Academic Excellence Project '5-100' and was also partially supported by the RFBR grant 19-01-00726-a.
 The work of D.T. was carried out within the framework of the State Programme of the Ministry of Education and Science of the Russian Federation, project 1.12873.2018/12.1, and was also partially supported by the RFBR grant 17-01-00366 A.

 \section{Compatible braidings}

Let $(R,F)$ be a couple of compatible braidings. As noticed in Introduction, the braiding $F$ is used for transferring the generating matrix $T$ (depending on parameters
or not) to the higher positions and the symmetry  $R$ comes in the defining relations of the algebras $\A(R,F)$ or that  $\Y(R,F)$.

We impose the following conditions on these braidings. We assume $R$ to be an involutive or Hecke  symmetry.
Remind  that a braiding $R$ is  called  an {\em involutive symmetry} (resp., a {\em  Hecke symmetry}) if it is
subject to the condition
  $$R^2=I\,\,\,({\rm resp.,}\,\,\, (R-q\, I)(R+\qq\, I)=0,\,\, q\in \C,\,\, q\not=\pm 1).$$

For such symmetries we construct current $R$-matrices according to the following {\em Baxterization procedure}.

\begin{proposition} (\cite{GS}) Let $R$ be an involutive or a Hecke symmetry. Define the operators
\be
R(u,v) = R -\frac{I}{u-v} \label{odin}
\ee
for an involutive symmetry $R$ and
\be
R(u,v) = R -\frac{(q-\qq)\, u\,  I}{u-v} \label{dva}
\ee
for a Hecke symmetry. The operators $R(u,v)$ are current $R$-matrices, i.e. they meet the Quantum Yang-Baxter equation with parameters
$$
R_{12}(u,v) R_{23}(u,w) R_{12}(v,w) =R_{23}(v,w) R_{12}(u,w) R_{23}(u,v).
$$
\end{proposition}

Namely, these current $R$-matrices are used in defining generalized Yangians. An $R$-matrix (\ref{odin}) (resp., (\ref{dva})) and the corresponding algebra $\Y(R,F)$
is called {\em rational} (resp., {\em trigonometric}).

As for the braidings $F$ we assume them to be {\em skew-invertible}. This means that there exists an operator $\Psi^F\in {\rm End}(V^{\otimes 2})$ such that
$$
{\rm Tr}_{(2)}F_{12}\Psi^F_{23} = {\rm Tr}_{(2)}\Psi^F_{12}F_{23} = P_{13}.
$$

If it is so,  we can introduce the so-called $F$-trace of any square $N\times N$ matrix $X$ by setting
$$
{\rm Tr}_FX = {\rm Tr}(C^F\cdot X), \qquad C^F = {\rm Tr}_{(2)}\Psi^F_{12}.
$$
Also, we put
$$
{\rm Tr}_{F(1\dots k)}X_{1\dots k} = {\rm Tr}_{(1\dots k)}(C_1^F\dots C_k^FX_{1\dots k})
$$
for any matrix $X_{1\dots k}$ of the appropriate size.

Now, analogically  to the notations $A_{\ov k}$  we introduce similar notations for
$N^2\times N^2$ matrices. Let  $A_{12}$ be such a matrix located at  positions number 1 and 2,  we put
$$A_{\ov {12}}=A_{{12}},\,\,     A_{\ov {13}}=F_{23} A_{\ov{12}} F_{23}^{-1},\,\, A_{\ov {23}}=F_{{12}}A_{\ov {13}}F_{{12}}^{-1}=
 F_{{12}}F_{{23}}A_{\ov{12}} F_{{23}}^{-1} F_{{12}}^{-1},     $$
and so on. In general, the notation $A_{\ov {kl}},\,\,k<l$ means that transferring of the matrix $A_{{12}}=A_{\ov{12}}$  to the positions number $k$ and $l$ is performed
by means of the  symmetry $F$ as follows
$$ A_{\ov {kl}}=(F_{k-1\, k}\,F_{k-2\, k-1} ...F_{12 }) (F_{l-1\, l}F_{l-2\, l-1}...F_{2 3})  A_{\ov {12}}
(F_{2 3}^{-1}... F_{l-2\, l-1}^{-1}\,F_{l-1\, l}^{-1})(F_{12}^{-1}  ...\,F_{k-2\, k-1}^{-1}  \,F_{k-1\, k}^{-1}). $$

Thus,  if  $(R,F)$ is a  couple of compatible braidings, the following holds
 \be   R_{23}=F_{12}\, F_{23}\, R_{12}\,F_{23}^{-1}\,F_{12}^{-1}=R_{\ov{23}}. \label{shift1} \ee
This means that  transferring of the braiding $R_{12}$  to the positions 2 and 3 performed either by means of  the usual flip $P$ or by means of the braiding
 $F$  leads to the same result: $R_{\ov{23}}=R_{23}$. This entails that $R_{\ov{i\, i+1}}=R_{i\, i+1}$ for any $i$.

It is clear  that the relation (\ref{braid}) for the operator $R$ is equivalent to the quantum Yang-Baxter equation
\be \RR_{12} \, \RR_{13}\, \RR_{23}= \RR_{23}\,\RR_{13}\, \RR_{12} \label{YB0} \ee
for the operator $\RR=R\,P$.

However, if  $(R,F)$ is a couple of compatible braidings,
 it is easy to see that the operator $\RR=R\, F$ is subject to the
following "braided version" of the  quantum Yang-Baxter equation
\be \RR_{\ov{12}} \, \RR_{\ov{13}}\, \RR_{\ov{23}}= \RR_{\ov{23}}\,\RR_{\ov{13}}\, \RR_{\ov{12}}. \label{YB} \ee

Moreover, the following generalization of (\ref{YB}) is valid
\be \RR_{\ov{ij}} \, \RR_{\ov{ik}}\, \RR_{\ov{jk}}= \RR_{\ov{jk}}\,\RR_{\ov{ik}}\, \RR_{\ov{ij}}, \label{YB1} \ee
provided $i,j,k$ are  positive integers such that $i<j<k$ .

The operators subject to (\ref{YB}) are called {\em braided $R$-matrices}.

Observe that if $F$ is an involutive symmetry,   the relation (\ref{YB1}) becomes valid for any positive pairwise distinct integers  $i,j,k$. Besides, the the notation $A_{\ov{ij}}$ is well-defined
for any matrix $A_{12}$ and any distinct positive integers $i$ and $j$. For instance, $A_{\ov{21}}=F\,A_{12}\,F$.

Consider a few examples of compatible braidings. The braidings  $R$ and $F=P$ are compatible for any $R$.  The corresponding algebra  $\Y(R,P)$ is called
the (generalized) Yangian of RTT type. The Drinfeld's Yangian $\Y(gl(N))$ is a particular case, respective to $R=P$. Another example of such an algebra is the so-called $q$-Yangian,
as it is defined in \cite{M}. The Hecke symmetry $R$ entering its definition is that coming from the quantum group $U_q(sl(N))$.

 It is also evident that if $F=R$, then the braidings $R$ and $F$ are compatible. We say that the corresponding braided  Yangian $\Y(R,R)$ is of Reflection Equation (RE) type.

If $F=P_{(m|n)}$ is a super-flip, and $R$ is the Hecke symmetry coming from the Quantum super-group $U_q(sl(m|n))$, the braidings $R$ and $F$ are compatible.
In the case $m=n=1$ the Hecke symmetry $R$ is represented in a basis by the following matrix
$$ \left(\begin{array}{cccc}
q&0&0&0\\
0&q-\qq&1&0\\
0&1&0&0\\
0&0&0&-\qq\end{array}\right).  $$

Note that this Hecke symmetry $R=R(q)$ is a deformation of the involutive symmetry $F=P_{(m|n)}$ and it depends analytically on $q$. Thus, we are in the frameworks of setting discussed at the end of Introduction.

\section{Symmetric polynomials  in generalized  Yangians}

In this section we define  (quantum) symmetric polynomials in generalized Yangians. Namely, we are dealing with elementary symmetric polynomials and power sums. However, first we consider
  an $N\times N$ numerical matrix $M$. In this case the  elementary symmetric polynomials $e_k(M)$ and quantum power sums $p_k(M)$ for this matrix
are respectively defined as follows
\be \det(M-t\, I)=\sum_0^N\,(-t)^{N-k}e_k(M),\,\,\, p_k(M)={\Tr}\, M^k. \label{classs} \ee
Note that if $M$ is triangular matrix, these  polynomials are respectively equal to the elementary symmetric polynomials and  power sums in eigenvalues of $M$.
This motivates the terminology.

If $M$ is a matrix with entries belonging to  a  noncommutative algebra,  it  is not in general possible to define analogs of these polynomials with interesting properties.
Fortunately, it is possible to do
in the algebras $\A(R,F)$ (see \cite{IOP})  and generalized Yangians $\Y(R,F)$. First, define elementary symmetric polynomials in the trigonometric generalized Yangians 

If $R$ is a Hecke symmetry, we put
\be
e_0(u)=1,\,\,\, e_k(u) = {\rm Tr}_{F(1\dots k)}A^{(k)}_{1\dots k}(R)\,T_{\overline 1}(u)T_{\overline 2}(q^{-2}u)\dots T_{\overline k}(q^{-2(k-1)}u),\,\,k=1,2...
\label{elem-symm}
\ee
Here, $A^{(k)}_{1\dots k}(R)$ is the skew-symmetrizer  acting in the space $V^{\ot k}$ and arising from the Hecke symmetry $R$. 
It can be defined by  the following recurrent relations
\be
A^{(1)}=I,\qquad
A^{(k+1)}_{1\dots k+1} = \frac{k_q}{(k+1)_q}\,A^{(k)}_{1\dots k}\Big(\frac{q^k}{k_q} \,I - R_k\Big)A^{(k)}_{1\dots k},\quad k\ge 1. \label{skew}
\ee

In the rational case  the corresponding elementary symmetric polynomials are defined as follows
\be
e_0(u)=1,\,\,\,   e_k(u) = {\rm Tr}_{F(1\dots k)}A^{(k)}_{1\dots k}(R)\,T_{\overline 1}(u)T_{\overline 2}(u-1)\dots T_{\overline k}(u-k+1),\,\,k=1,2...
\label{elem-symm-inv}
\ee
Here $A^{(k)}_{1\dots k}(R)$ is the skew-symmetrizer respective to the involutive symmetry $R$. Its explicit formula can be obtained from that (\ref{skew}) at $q=1$.

As for the power sums, we define them respectively as follows
$$
p_0(u)=1,\quad  p_k(u) = {\rm Tr}_{F(1\dots k)}T_{\overline 1}(q^{-2(k-1)}u)T_{\overline 2}(q^{-2(k-2)}u)\dots T_{\overline k}(u)R_{k-1}\dots R_2R_1,\quad k\ge 1,
$$
$$
p_0(u)=1,\quad p_k(u) = {\rm Tr}_{F(1\dots k)}T_{\overline 1}(u-k+1)T_{\overline 2}(u-k+2)\dots T_{\overline k}(u)R_{k-1}\dots R_2R_1,\quad k\ge 1.
$$

It should be emphasized that in the braided Yangians of RE type these formulae could be reduced to the following forms respectively
$$
p_k(u) = {\rm Tr}_{F} T(q^{-2(k-1)}u)T(q^{-2(k-2)}u)\dots T(u),\qquad
p_k(u) = {\rm Tr}_{F} T(u-k+1)T(u-k+2)\dots T(u).
$$
Observe that these formulae are in a sense similar to the second formula from (\ref{classs}) but they contain shifts of the arguments, multiplicative and additive respectively.

Also, note that the elementary symmetric polynomials and power sums are related via a quantum  version of the Newton identities. If $R$ is a Hecke symmetry these identities are
$$
k_qe_k(u) - q^{k-1}p_1(q^{-2(k-1)}u)e_{k-1}(u)+q^{k-2}p_2(q^{-2(k-2)}u)e_{k-2}(u)+\dots +(-1)^kp_k(u)e_{0}(u)= 0.
\label{q-Newt}
$$
If  $R$ is involutive,  then we have
$$
k e_k(u) - p_1(u-k+1)e_{k-1}(u)+p_2(u-k+2)e_{k-2}(u)+\dots +(-1)^kp_k(u)= 0.
$$

A proof of these identities is given in \cite{GS} for the braided Yangians of RE type.  For braided Yangians of general form these identities can be shown in a similar way.

Also, note that if the bi-rank of $R$ is $(m|0)$, then $e_k(u)\equiv 0$ for $k>m$.

The subalgebra generated in $\Y(R,F)$ by the elements $e_k(u)$ is called {\em Bethe} one.

\begin{proposition}
Let $\Y(R,F)$ be a rational or trigonometric generalized Yangian. Then the elements $e_k(u), \quad k=1,2... $ and consequently these $p_k(u),\quad k=1,2... $
commute with each other:
$$e_k(u)\,e_p(v)= e_p(v)\,e_k(u), \qquad \forall \, k,p,\quad \forall\,u,v, $$
and consequently the Bethe subalgebra is commutative.
\end{proposition}

 A detailed proof of this claim is given  in \cite{GSS} for the braided Yangians of RE type. Generalized  Yangians corresponding to all couples $(R,F)$ under consideration can be treated in a similar way.

\begin{remark} {\rm In  \cite{IO} a notion of half-quantum algebras (HQA) was introduced. Each of
these algebras is
also defined via a couple  of compatible braidings,  namely, by the system
\be
S^{(2)}_{12}(R)\, T_{\ov 1}\, T_{\ov 2}\,A^{(2)}_{12}(R)=0,\,\,  {\rm where}\,\, A^{(2)}_{12}(R)=\frac{q\, I-R}{q+\qq},\,\, S^{(2)}_{12}(R)=\frac{\qq\, I+R}{q+\qq}
\label{half}
\ee
are the skew-symmetrizer  and  symmetrizer respectively, provided $R$ is a Hecke symmetry. As usual, the braiding $F$ is employed  for defining the overlined indexes.

Let us exhibit an equivalent form of (\ref{half}), which is useful in the study of the generalized Yangians,
$$
A^{(2)}_{12}(R)\, T_{\ov 1}\, T_{\ov 2}\,A^{(2)}_{12}(R)=T_{\ov 1}\, T_{\ov 2}\, A^{(2)}_{12}(R). $$

In the HQA there exist analogs of the elementary symmetric polynomials and power sums and those of the Newton and Cayley-Hamilton identities (see \cite{IO}).

However, in general in the HQA the commutativity of the these symmetric polynomials is not valid.

The HQA are related to the braided Yangians as follows.
In the  trigonometric $R$-matrix $R(u,v)$
we put $v=q^{-2}\, u$. Then  we have
$$ R(u, q^{-2} u)= R(q)-\frac{(q-\qq) q^2\, I}{(q^2-1)}=R(q)-q\, I. $$
This operator  coincides  up to a factor $-(q+\qq)$ with the skew-symmetrizer $A^{(2)}(R)$.

After multiplying the defining system of the corresponding braided Yangian evaluated at  $v=q^{-2}\, u$
by $S^{(2)}_R$ from the right hand side, we get the relation
$$A^{(2)}_{12}(R)\, T_{\ov 1}(u)\, T_{\ov 2}(q^{-2} u)\,S^{(2)}_{12}(R)=0, $$
which can be written under the following form
$$A^{(2)}_{12}(R)\, (q^{-2 u\pa_u}\,T_{\ov 1}(u))\,(q^{-2 u\pa_u}\,T_{\ov 2}(u))\,S^{(2)}_{12}(R)=0, \,\, \pa_u=\frac{d}{d\, u}.$$
This relation looks like that in a HQA but the role of the matrices $T(u)$ is played  by the operator $q^{-2 u\pa_u}\,T(u)$.

A similar treatment is possible in the rational braided Yangians but in them the shifts are additive, since the parameters $u$ and $v$ are related as follows $u-v=1$.  }

\end{remark}

\section{Braided $r$-matrices and braided Sklyanin brackets}

Let $F$ be an involutive symmetry.  Let us consider the following operator
 \be r(u,v)=\frac{F}{u-v}, \label{brmat} \ee
 which is a braided generalization of the rational $r$-matrix $r(u,v)=\frac{P}{u-v}$.

 It is easy to see that it meets  the following relations
  \be 1.\,\,  r_{\ov{21}}(v,u)+r_{\ov{12}}(u,v)=0, \label{fir}  \ee
 \be 2.\,\, [r_{\ov {12}}(u,v)\, r_{\ov {13}}(u,w)]+[r_{\ov {12}}(u,v)\, r_{\ov {23}}(v,w)]+[r_{\ov {13}}(u,w)\, r_{\ov {23}}(v,w)]=0.\label{sec} \ee

In order to define a braided analog of   {\em trigonometric} $r$-matrix we need two compatible braidings.
Again, let $(R,F)$ be a  couple of compatible braidings. Also,  suppose  that $F$ is an involutive symmetry and $R=R(q)$ is a Hecke symmetry deforming $F$ as
 described in Introduction. Let us  expand the operator $\RR=R\, F$ at the point $q=1$:
\be  \RR= I+h\, r+O(h^2),\, q=exp(h). \label{rm} \ee

\begin{definition} The element $r\in End(\vv)$  entering
 this expansion is called a (constant) braided $r$-matrix. \end{definition}

\begin{proposition} This braided $r$-matrix $r$ has the following properties
\be 1. \,\, r_{\ov {12}}+r_{\ov{21}}=2F,  \label{firr}  \ee
\be 2. \,\, [r_{\ov {12}}\, r_{\ov {13}}]+[r_{\ov {12}}\,r_{\ov {23}}]+[r_{\ov {13}}\, r_{\ov {23}}]=0.  \label{secc}\ee
\end{proposition}

{\bf Proof} Similarly to the classical case, the relation (\ref{secc}) immediately follows from (\ref{YB}).  The relation (\ref{secc}) follows from that
$$ (\RR\, F)^2=(q-\qq)(\RR\, F)+I.$$

Now,   consider the following  operator
 \be r(u,v)=\frac{F\,u}{ u-v}-\frac{r}{2}. \label{brmat1} \ee
We call this operator  {\em braided trigonometric $r$-matrix}.  As usual, the term "trigonometric" is justified by  another form of this operator
obtained by the change $u\to q^u,\,\, v\to q^v$.

\begin{proposition} \label{four1} The operator (\ref{brmat1})  meets the  relations (\ref{fir}) and (\ref{sec})
\end{proposition}

{\bf Proof} The first relation follows immediately from (\ref{firr}). In order to show (\ref{sec}),
we have to compute the braided Schouten bracket of the operator $r(u,v)$ with itself.
Let us precise that by the {\em braided Schouten bracket} of two such operators $A(u,v)$ and $B(u,v)$ we mean the following expression
$$[[A,B]](u,v,w)=[A_{\ov{12}}(u,v), B_{\ov{13}}(u,w)]+[A_{\ov{12}}(u,v), B_{\ov{23}}(v,w)]+[A_{\ov{13}}(u,w), B_{\ov{23}}(v,w)]+  $$ $$
[B_{\ov{12}}(u,v), A_{\ov{13}}(u,w)]+[B_{\ov{12}}(u,v), A_{\ov{23}}(v,w)]+[B_{\ov{13}}(u,w), A_{\ov{23}}(v,w)].$$
If $A$ and/or $B$ are constant, the corresponding parameters should be omitted.

By direct computations we have that the Schouten bracket of the summand $A=\frac{F\, u}{u-v}$ with itself is equal to
$$[[A,A]](u,v,w)=\frac{2\,u}{u-w} [F_{23},F_{12}]. $$
Also, the following holds
$$[[A, r]](u,v,w)=[\frac{F_{\ov{ 12}}\, u}{u-v}, r_{\ov{13}}]+[r_{\ov{12}}, \frac{F_{\ov{13}}\, u}{u-w}]+[\frac{F_{\ov{12}}\, u}{u-v}, r_{\ov{23}}]+
[r_{\ov{12}}, \frac{F_{\ov{23}}\, v}{v-w}]+ [\frac{F_{\ov{ 13}}\, u}{u-w}, r_{\ov{23}}]+[r_{\ov{13}}, \frac{F_{\ov{23}}\, v}{v-w}]$$
$$=[F_{\ov{ 13}},r_{\ov{23}}-r_{\ov{12}}]\frac{u}{u-w}=(-r_{\ov{23}}+r_{\ov{12}}+r_{\ov{21}}-r_{\ov{32}})\,\frac{\, F_{\ov{13}}u}{u-w}=2[F_{23}, F_{12}]\frac{u}{u-w}.   $$
Besides, $[[r, r]]=0$ in virtue of (\ref{secc}). This completes the proof.

Let us remark that due to the  property (\ref{fir}) of the braided $r$-matrix $r(u,v)$ defined by  (\ref{brmat1}) it can be  cast under the following form
$$r(u,v)=\frac{1}{2}\left(\frac{F(u+v)}{u-v}-\frac{r-r_{\ov{21}}}{2}\right).$$

It should be emphasized that we do not use any concrete form of the symmetries $R$ and $F$.

Below, we also need following properties of braided $R-$ and $r$-matrices.

\begin{proposition} If $(R,F)$ is a couple of compatible braidings, then the following holds
$$ \RR_{12} \, A_{\ov 3}=A_{\ov 3}\, \RR_{12}$$
for any $N\times N$ matrix $A$.  \end{proposition}
 {\bf Proof} By using the compatibility of the braidings $R$ and $F$, we get
$$\RR_{12}\, F_{23} \, F_{12}\, A_1\, F_{12}^{-1}\, F_{23}^{-1}=R_{12}\, F_{12}\,F_{23} \, F_{12}\, A_1\, F_{12}^{-1}\, F_{23}^{-1}=
R_{12}\, F_{23}\,F_{12} \, F_{23}\, A_1\, F_{12}^{-1}\, F_{23}^{-1}=$$
$$F_{23}\,F_{12}\,R_{23}\,F_{23}\,A_1\, F_{12}^{-1}\, F_{23}^{-1}=
F_{23}\,F_{12}\,A_1\,R_{23} F_{23}\, F_{12}^{-1}\, F_{23}^{-1}=F_{23} \, F_{12}\, A_1\, F_{12}^{-1}\, F_{23}^{-1}\, \RR_{12}.$$

\begin{corollary} Additionally,  suppose  $F$ to be an involutive symmetry.  Then for any pairwise distinct positive integers $i,j,k$ we have
\be \RR_{\ov{ij}} \, A_{\ov k}=A_{\ov k}\, \RR_{\ov{ij}}.
\label{commutt} \ee
\end{corollary}

\begin{proposition} Under the same hypothesis, if   $A$ is an  $N\times N$ matrix and $r(u,v)$ is a braided current (rational or trigonometric) $r$-matrix, then the following holds
\be r_{\ov{ij}}(u,v)\, A_{\ov k} = A_{\ov k}\, r_{\ov{ij}}(u,v). \label{comm} \ee
\end{proposition}

{\bf Proof} It suffices to consider the case $i=1,\, j=2,\, k=3$. Then this relation is clear for any braided rational $r$-matrix. It is also so for the first summand of any braided trigonometric $r$-matrix. For the second summand it follows from the previous claim.

Now, we pass to describing a braided analog of the first Sklyanin bracket.

Let $F$ be in involutive symmetry and
 $r(u,v)$ be a braided rational or trigonometric  $r$-matrix defined correspondingly by (\ref{brmat}) or (\ref{brmat1}).

Let us define the following involutive symmetry on the space of the currents as follows
$$\FF(T_{\ov 1}(u)\ot T_{\ov 2}(v))=F\,(T_{\ov 1}(v)\ot T_{\ov 2}(u))\, F=T_{\ov 2}(v)\ot T_{\ov 1}(u).$$

Now, introduce the following Lie type operator, which is a braided analog of the first Sklyanin bracket
\be  [T_{\ov 1}(u), T_{\ov 2}(v)]^{\FF}=[T_{\ov 1}(u)+ T_{\ov 2}(v), r(u,v)].
\label{bbrr} \ee

\begin{proposition} The operator (\ref{bbrr}) meets the following conditions

1. $[T_{\ov 1}(u), T_{\ov 2}(v)]^{\FF}=-[\,,\,]^{\FF} \FF(T_{\ov 1}(u), T_{\ov 2}(v))=-[\,,\,]^{\FF} (T_{\ov 2}(v), T_{\ov 1}(u)),$

2. $[\,,\,]^{\FF}\,[\,,\,]_{23}^{\FF}(I+\FF_{12}\,\FF_{23}+\FF_{23}\,\FF_{12})(T_{\ov 1}(u)\ot T_{\ov 2}(v)\ot T_{\ov 3}(w))=0.$
\end{proposition}

{\bf Proof} The first claim  follows immediately from the skew-symmetry of the braided $r$-matrix $r(u,v)$. In order
to prove the second one we use the  relations (\ref{comm}) with $A=T$. Then we have
$$[T_{\ov 1}(u),  [T_{\ov 2}(v),  T_{\ov 3}(w)]^\FF]^\FF=[T_{\ov 1}(u),[T_{\ov 2}(v)+ T_{\ov 3}(w), r_{\ov{23}}(v,w)]]^\FF=$$
$$[T_{\ov 1}(u), (T_{\ov 2}(v)+ T_{\ov 3}(w))\, r_{\ov{23}}(v,w)- r_{\ov{23}}(v,w)\, (T_{\ov 2}(v)+ T_{\ov 3}(w))]^\FF=$$
$$([T_{\ov 1}(u),T_{\ov 2}(v)]^\FF+ [T_{\ov 1}(u),T_{\ov 3}(w)]^\FF)\, r_{\ov{23}}(v,w)-r_{\ov{23}}(v,w)\,([T_{\ov 1}(u), T_{\ov 2}(v)]^\FF+
[T_{\ov 1}(u),T_{\ov 3}(w)]^\FF).$$
By using the relation (\ref{comm})  once more we  can reduce this expression to the following form
 $$[\XX, r_{\ov{12}}(u,v)]\, r_{\ov{23}}(v,w)+ [\XX, r_{\ov{13}}(u,w)]\, r_{\ov{23}}(v,w)-r_{\ov{23}}(v,w)\,
[\XX, r_{\ov{12}}(u,v)]+ r_{\ov{23}}(v,w) [\XX, r_{\ov{13}}(u,w)],$$
where $\XX=T_{\ov 1}(u)+T_{\ov 2}(v)+ T_{\ov 3}(w)$.

Now, by applying the operators $\FF_{\ov{12}}\, \FF_{\ov{23}}$ and $ \FF_{\ov{23}}\,\FF_{\ov{12}}$ to this expression and by adding all results we arrive to the conclusion.

\section{Braided   KZ equations--classical and quantum}

Let us pass to constructing a family of commuting differential operators looking like the famous KZ ones.  To this end we consider two families of matrices
 $$M_i= g_{\ov i}+\kappa \sum_{k\not=i}^n \frac{F_{\ov{i\, k}}}{u_i-u_k},\,\, i=1,...,n$$
 and
  $$N_i=g_{\ov i}+\kappa \sum_{k\not=i}^n (\frac{F_{\ov{i\, k}}\, u_i}{u_i-u_k}-\frac{r_{\ov{i\, k}}}{2}),\,\, i=1,...,n,$$
  associated with  braided rational and trigonometric $r$-matrices correspondingly.
Here   $n\geq 2$ is an integer, $\kappa\in\C$ is an arbitrary parameter,  $u_i\in \C,i=1,...,n$ are pairwise distinct numbers, and $g$ is a numerical  $N\times N$ matrix.
Also, throughout this section $F$ is assumed to be a skew-invertible involutive symmetry.

Also, constitute two families of differential operators
\be \pa_i-M_i=\pa_i-\left(g_{\ov i}+\kappa \sum_{k\not=i}^n \frac{F_{\ov{i\, k}}}{u_i-u_k}\right),\,\, i=1,2,...,n, \label{KZ1} \ee
and
\be u_i\, \pa_i-N_i=u_i\, \pa_i-\left(g_{\ov i}+\kappa \sum_{k\not=i}^n (\frac{F_{\ov{i\, k}}\, u_i}{u_i-u_k}-\frac{r_{\ov{i\, k}}}{2})\right),\,\, i=1,2,...,n. \label{KZ2} \ee
where $\pa_i=\pa_{u_i}$.

\begin{proposition}
Let us assume that the matrix $g$ is subject to the relation $g_{\ov 1}\, g_{\ov 2}=g_{\ov 2}\, g_{\ov 1}$. Then the following holds true
\be  \pa_i\, M_j-\pa_j\, M_i-[M_i, M_j]=0,\,\,\,\,\, u_i\pa_i\, N_j-u_j\pa_j\, N_i-[N_i, N_j]=0 . \label{pro} \ee
\end{proposition}

We call the operators (\ref{KZ1}) and (\ref{KZ2}) braided KZ connections--rational and trigonometric respectively.

First, observe that the condition on the matrix $g$ can be cast under the following form
\be F_{12}\, g_1\, F_{12}\, g_1=g_1\,F_{12}\, g_1\, F_{12}, \label{gg} \ee
which means that this matrix realizes a one-dimensional representation of the RE algebra associated with  the involutive symmetry $F$.

{\bf Proof} of the first relation from (\ref{pro}) results  from the fact that the operator $\frac{F}{u-v}$ is a braided $r$-matrix and the
following relations
\be  \pa_i \left(\frac{F_{\ov{j\, i}}}{u_j-u_i}\right)=\pa_j \left(\frac{F_{\ov{i\, j}}}{u_i-u_j}\right) ,\label{KZ3} \ee
$$[g_{\ov i}, F_{\ov{j\, k}}]=0,\,\,  [g_{\ov i}, F_{\ov{j\, i}}]+[g_{\ov j}, F_{\ov{i\, j}}]=0,\,\, [g_{\ov i}, g_{\ov j}]=0,$$
where $i,\, j,\, k$ are pairwise distinct. (Note that $F_{\ov{i\, j}}=F_{\ov{j\, i}}$.)

The second relation  from (\ref{pro}) can be proven in the same way,  with using the relation (\ref{comm}) and the following one
$$ u_i\pa_i \left(\frac{F_{\ov{j\, i}}u_j}{u_j-u_i}\right)=u_j\pa_j \left(\frac{F_{\ov{i\, j}}u_i }{u_i-u_j}\right) $$
instead of (\ref{KZ3}).

\begin{corollary} The corresponding systems of differential equations (called the KZ ones)
$$\pa_i\, \Psi=M_i \,\Psi,  \,\, i=1,2...n,$$
$$u_i\,\pa_i\, \Psi=N_i \,\Psi, \,\, i=1,2...n,$$
where $\Psi$ is a vector-function of the length $N$,
are compatible.
\end{corollary}

Our next purpose is to introduce quantum counterparts of these systems. In the case related to the affine Quantum Groups these counterparts were introduced in \cite{FR}. Similarly to
\cite{FR} (also, see \cite{EFK}), we get systems of difference equations. However, we deal without  any Quantum Group. Instead, we use the  rational and trigonometric braidings $R(u,v)$
defined by (\ref{brmat}) or (\ref{brmat1}) correspondingly.

Consider the following operators
\be  \RR(u,v)=R(u,v)\, F\, f(u,v)^{-1} \label{RR} \ee
where
\be  f(u,v)=1-\frac{1}{u-v} \,\,\, {\rm{or}}\,\,\, f(u,v)=q-\frac{(q-\qq)\, u}{u-v} \label{fuv} \ee
The former (resp., latter) $f(u,v)$ corresponds to  a rational (resp., trigonometric) case.

Observe that due to the factor $f(u,v)^{-1}$ in (\ref{RR})  the following relation
$$\RR_{\ov{12}}(u,v)\, \RR_{\ov{21}}(v,u)=I$$
holds.

Similarly to the classical braided KZ operators above, below we use a matrix $g$ subject to the  condition $g_{\ov i}\, g_{\ov j}=g_{\ov j}\,g_{\ov i}$. Observe that
in virtue of (\ref{commutt}) we have
$$[R_{\ov{ij}}, g_{\ov k}]=0,$$
Here $i,j,k$ are positive integer pairwise distinct.
Let us additionally demand
$$[R_{\ov{ij}}, g_{\ov i}\ot g_{\ov j}]=0.$$

Besides, introduce  operators
\be T_i^{add} f(u_1...u_i...u_n)=f(u_1...u_i+p...u_n),\,\,T_i^{mult} f(u_1...u_i...u_n)=f(u_1...p\,u_i...u_n) \label{diff} \ee
which perform the shifts of the variables, additive and multiplicative respectively. Here, $p\in \C$ is an arbitrary nontrivial number.

Observe that any rational (resp., trigonometric) $R$-matrix  $R(u,v)$ is invariant with respect to the operator $T_i^{add}$ (resp., $T_i^{mult}$) applied to the both parameters
$u$ and $v$. Namely, we have
$$(T_i\ot T_i)(R(u,v))=R(u,v).$$
Hereafter,  $T_i$ for $i=1,2$ stands for the operator $T_i^{add}$ or that $T_i=T_i^{mult}$ in function of the type of $R(u,v)$.

Let us introduce the following notations
$$\RR_{i\, \downarrow}= \RR_{\ov{i\, i-1}}(u_i,u_{i-1})...\RR_{\ov{i\,1}}(u_i,u_{1}),$$
$$\RR_{i\, \uparrow}=\RR_{\ov{i\, n}}(u_i,u_{n})...\RR_{\ov{i\, i+1}}(u_i,u_{i+1}).$$
Also, introduce the following operators
\be
\label{delta}
\Delta_i=\RR_{i\, \downarrow}\, T_i^{-1}\, g_{\ov i}\, \RR_{i\, \uparrow}=\RR_{i\, \downarrow}\, \Theta_i \, \RR_{i\, \uparrow}.
\ee
Hereafter, we use the following notation $\Theta_i=T_i^{-1}\, g_{\ov i}$. 

We are interested in the compatibility condition of the system
\be T_i \, \Psi(u_1...u_n)= \kappa\, T_i\, \RR_{i\, \downarrow}\, \Theta_i\, \RR_{i\, \uparrow} \, \Psi(u_1...u_n),\,\, i=1...n, \label{QQ}\ee
where $\kappa$ is an arbitrary nontrivial number.

More explicitly, in the rational case the $i$-th equation of this system reads
\be
\label{system}
\Psi(u_1...u_i+p...u_n) =\kappa\mathcal{M}_i \Psi(u_1...u_i...u_n)
\ee
where
$$\mathcal{M}_i=T_i\,\RR_{i\, \downarrow}\,\Theta_i\,\RR_{i\, \uparrow} =\RR_{\ov{i\, i-1}}(u_i+p,u_{i-1})...\RR_{\ov{i\,1}}(u_i+p,u_{1})g_{\ov i}\RR_{\ov{i\, n}}(u_i,u_{n})... \RR_{\ov{i\, i+1}}(u_i,u_{i+1}).$$
The compatibility condition of the system (\ref{system}) is called the holonomy condition (see \cite{EFK}, section 10.5). It takes the form
\be
\label{commm}
T_i \mathcal{M}_j T_i^{-1} \mathcal{M}_i = T_j \mathcal{M}_i T_j^{-1} \mathcal{M}_j.
\ee
The commutativity of $T_i$ and $T_j$ provides an equivalent form of (\ref{commm})
\bea
[ T_j^{-1} \mathcal{M}_j, T_i^{-1} \mathcal{M}_i]=0.\nn
\eea

\begin{proposition} The system (\ref{QQ}) satisfies the holonomy condition in both rational and trigo\-no\-met\-ri\-c cases.
\end{proposition}
{\bf Proof.}
The demonstration is similar to the non-braided case. We exhibit it to make our paper self-contained. 

Below, we systematically use the following relations
$$\RR_{\ov {ij}}\, \RR_{\ov {kl}}=\RR_{\ov {kl}}\,\RR_{\ov {ij}},$$
provided $i,j,k,l$ are positive pairwise distinct integers. 

\begin{figure}[h!]
\center
\includegraphics[width=80mm]{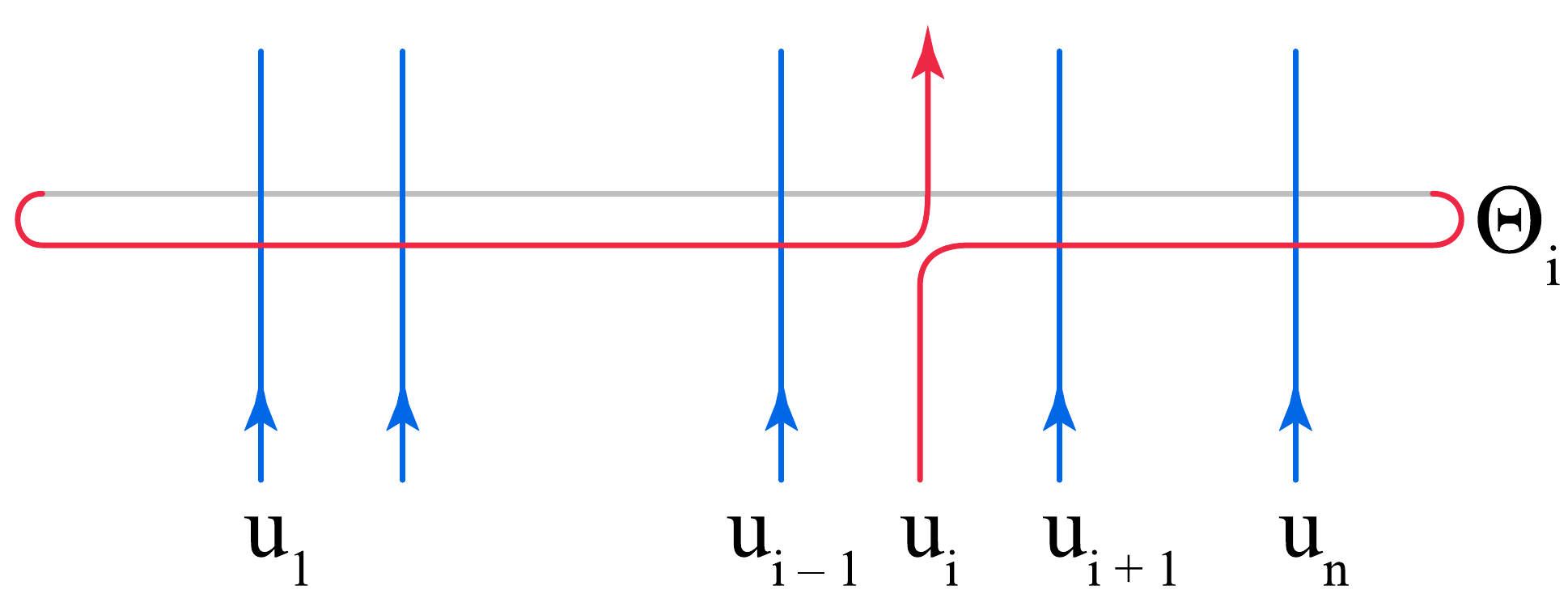}
\caption{$\Delta_i$}
\label{pic-deltai}
\end{figure}

\begin{figure}[h!]
\center
\includegraphics[width=40mm]{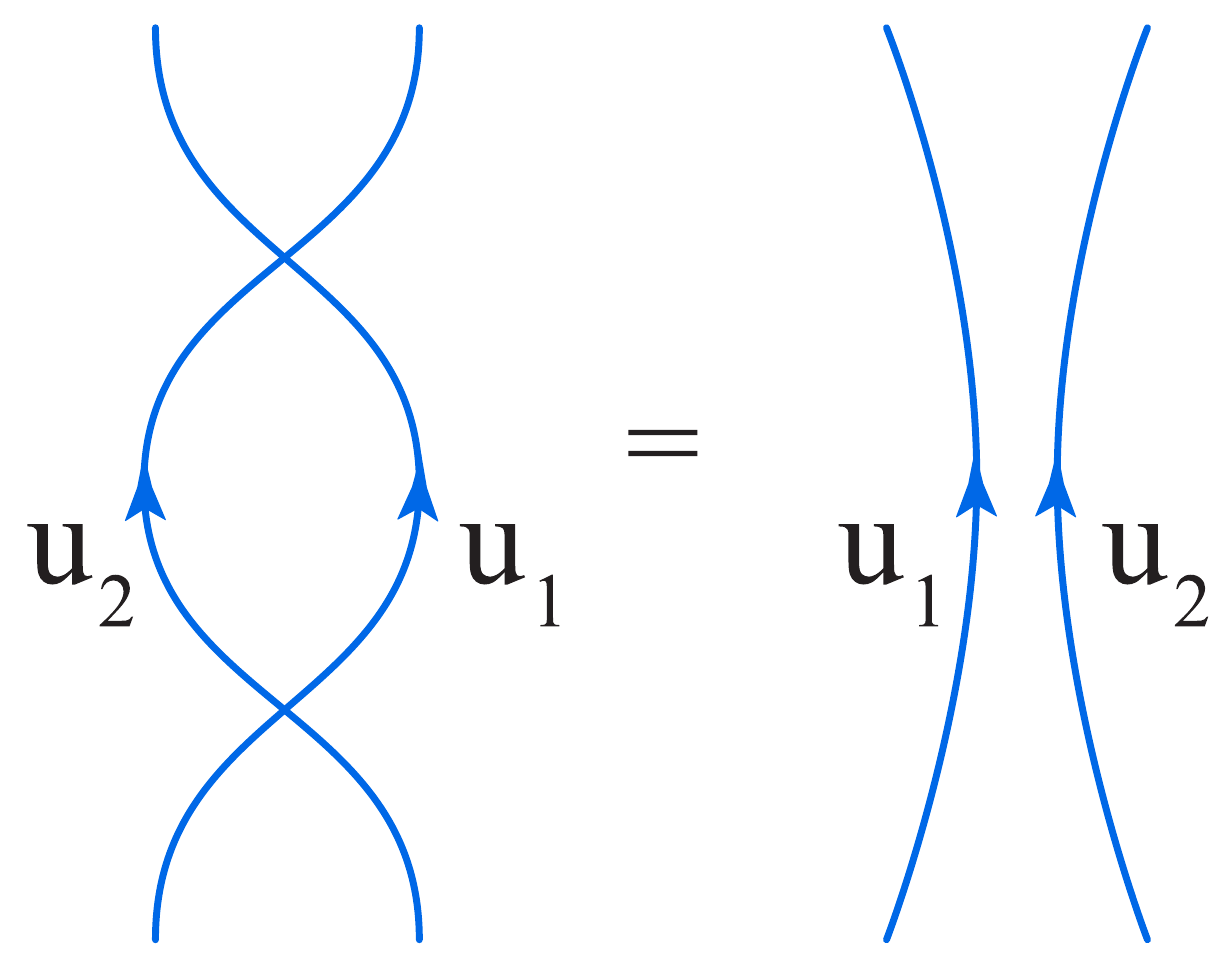}
\caption{2-nd Reidemeister move}
\label{pic-R2}
\end{figure}

\begin{figure}[h!]
\center
\includegraphics[width=50mm]{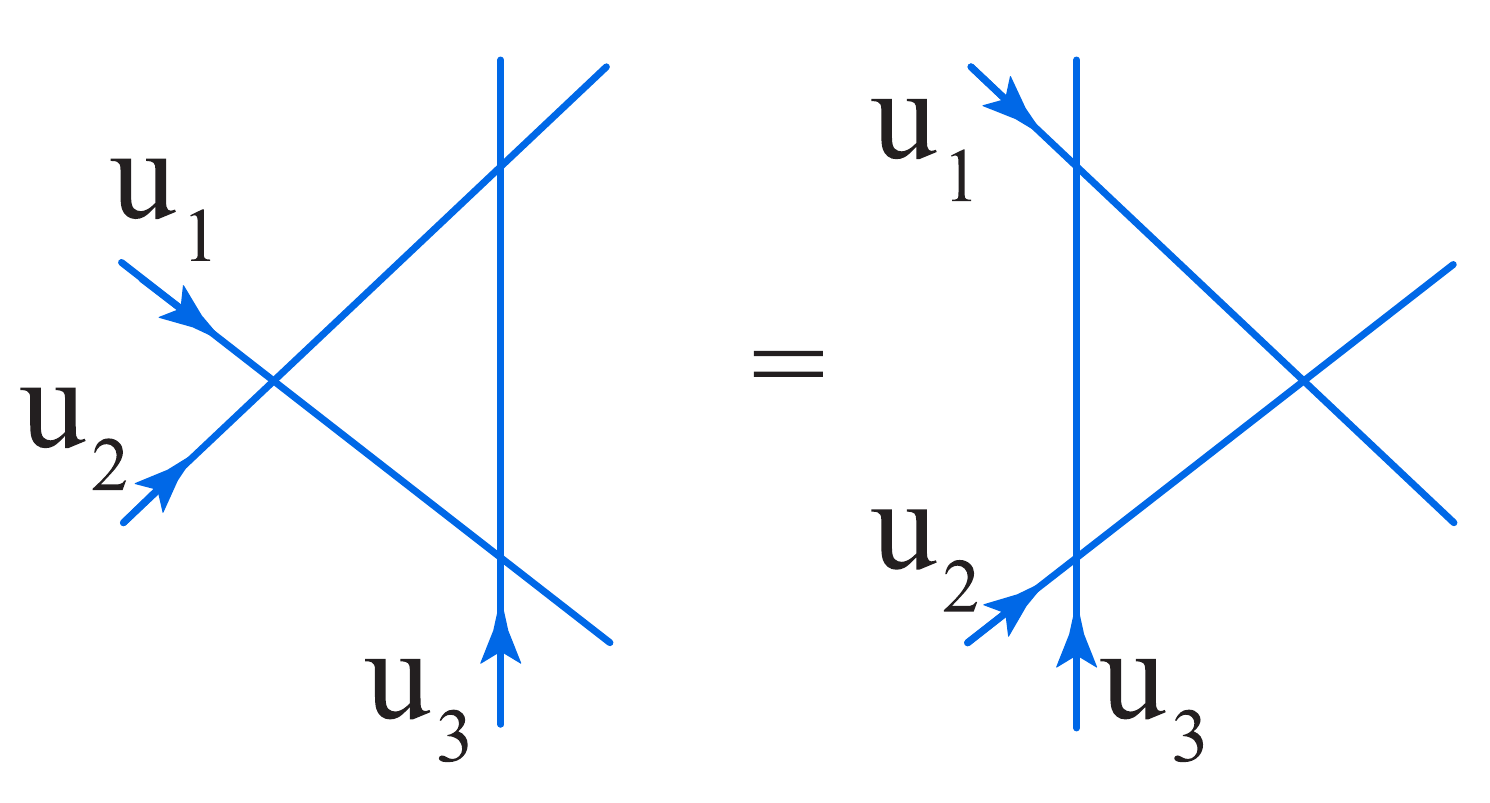}
\caption{YBE}
\label{pic-R3}
\end{figure}

Let us illustrate the expression $\Delta_i$ by the figure \ref{pic-deltai}.
Here, we use the standard pictorial representation for the basic structure equations involving $\RR$-operator. These resemble the 2-nd and 3-rd Reidemeister moves illustrated in figures \ref{pic-R2} and \ref{pic-R3}.

All pictures imply the application of $\RR$ at the crossing points to the vector spaces with  corresponding numbers and with appropriate spectral parameters marked on the diagram. Here the symbol  $\Theta_i=T_i^{-1} g_{\ov i}$ means the application of this operator to the corresponding space.

\begin{figure}[h!]
\center
\includegraphics[width=100mm]{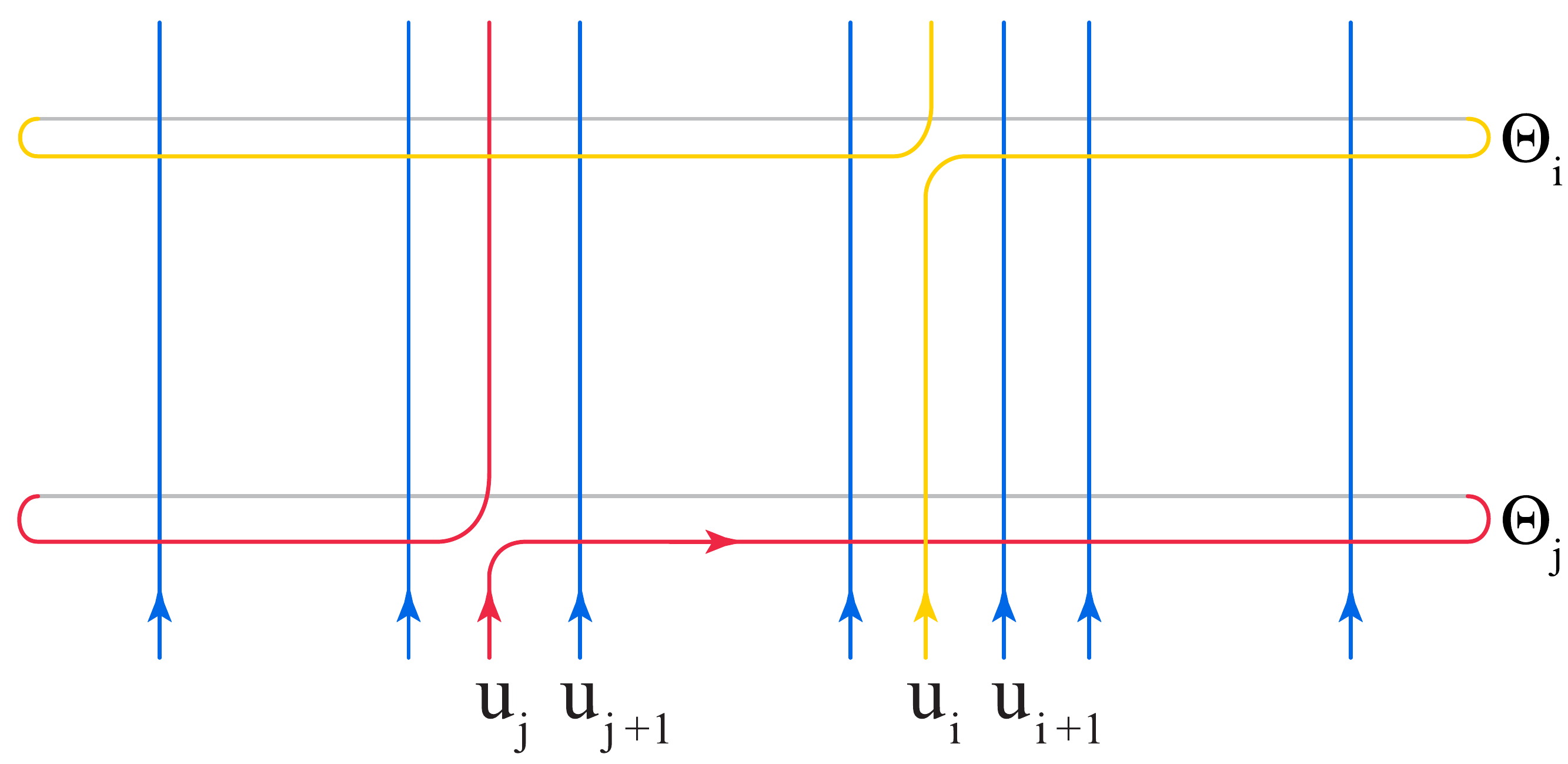}
\caption{${\Delta_i} \Delta_j$}
\label{pic-deltaij}
\end{figure}

On the figure \ref{pic-deltaij} we draw the diagram, corresponding to the expression 
${\Delta_i} \Delta_j.$
Without loss of generality we assume that $j<i.$ We perform a series of mutations using the Yang-Baxter equation (YBE) and once the commutation relation with the shift operator.
\begin{figure}[h!]
\center
\includegraphics[width=100mm]{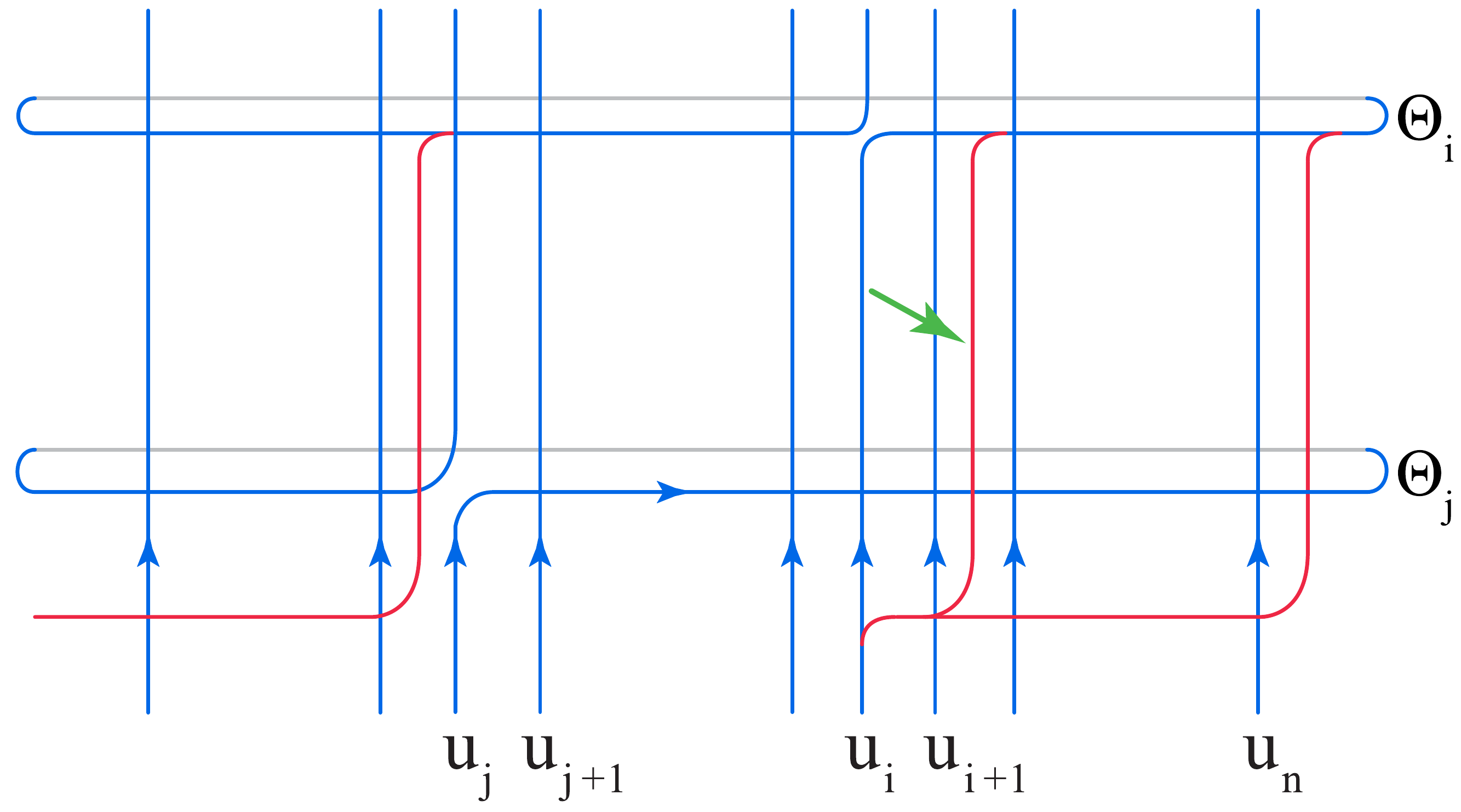}
\caption{}
\label{pic-deltaijmut1}
\end{figure}
The green arrow on the figure \ref{pic-deltaijmut1} illustrates the following algebraic transformation. First, we push forward the element ${\color{red}\RR_{\ov{i\, i+1}}}$ till ${\color{blue}\RR_{\ov{j\, i+1}}\RR_{\ov{j\, i}} }.$
\bea
\Delta_i \Delta_j&=&\RR_{\ov{i\, i-1}}\ldots \RR_{\ov{i\, 1}}\Theta_i \RR_{\ov{i\,n}}\ldots {\color{red}\RR_{\ov{i\, i+1}}}\RR_{\ov{j\, j-1}}\ldots \RR_{\ov{j\, 1}} \Theta_j\RR_{\ov{j\,n}}\ldots \RR_{\ov{j\,j+1}}\nn\\
&=&\RR_{\ov{i\, i-1}}\ldots \RR_{\ov{i\, 1}}\Theta_i\RR_{\ov{i\,n}}\ldots \RR_{\ov{i\, i+2}}\RR_{\ov{j\, j-1}}\ldots \RR_{\ov{j\, 1}} \Theta_j\RR_{\ov{j\,n}}\ldots\RR_{\ov{j\, i+2}}\nn\\
&\times&{\color{red}\RR_{\ov{i\, i+1}}}{\color{blue}\RR_{\ov{j\, i+1}}\RR_{\ov{j\, i}} }\ldots\RR_{\ov{j\,j+1}}.\nn
\eea
Now, by using the YBE we arrive to the following formula
\bea
\Delta_i \Delta_j&=&\RR_{\ov{i\, i-1}}\ldots \RR_{\ov{i\, 1}}\Theta_i\RR_{\ov{i\,n}}\ldots \RR_{\ov{i\, i+2}}\RR_{\ov{j\, j-1}}\ldots \RR_{\ov{j\, 1}}\Theta_j\RR_{\ov{j\,n}}\ldots\RR_{\ov{j\, i+2}}\nn\\
&\times&{\color{blue}\RR_{\ov{j\, i}}\RR_{\ov{j\, i+1}}}\RR_{\ov{j\,i-1}}\ldots\RR_{\ov{j\,j+1}}{\color{red}\RR_{\ov{i\, i+1}}}.\nn
\eea

Then, by the similar transformation we transpos the groups of magenta elements and the group of blue ones.
\bea
\Delta_i \Delta_j&=&\RR_{\ov{i\, i-1}}\ldots \RR_{\ov{i\, 1}} \Theta_i{\color{magenta}\RR_{\ov{i\,n}}\ldots \RR_{\ov{i\, i+2}}}\RR_{\ov{j\, j-1}}\ldots \RR_{\ov{j\, 1}} \Theta_j{\color{blue}\RR_{\ov{j\,n}}\ldots\RR_{\ov{j\, i+2}}}\nn\\
&\times&{\color{green}\RR_{\ov{j\, i}}}\RR_{\ov{j\, i+1}} \RR_{\ov{j\,i-1}}\ldots\RR_{\ov{j\,j+1}}\RR_{\ov{i\, i+1}}\nn\\
&=&\RR_{\ov{i\, i-1}}\ldots \RR_{\ov{i\, 1}} \Theta_i\RR_{\ov{j\, j-1}}\ldots \RR_{\ov{j\, 1}} \Theta_j{\color{green}\RR_{\ov{j\, i}}}{\color{blue}\RR_{\ov{j\,n}}\ldots\RR_{\ov{j\, i+2}}\RR_{\ov{j\, i+1}} }\nn\\
&\times&
\RR_{\ov{j\,i-1}}\ldots\RR_{\ov{j\,j+1}}{\color{magenta}\RR_{\ov{i\,n}}\ldots \RR_{\ov{i\, i+2}}}\RR_{\ov{i\, i+1}}.\nn
\eea
In this procidure we used the YBE involving the green operator ${\color{green}\RR_{\ov{j\,i}}}$ several times.
\begin{figure}[h!]
\center
\includegraphics[width=100mm]{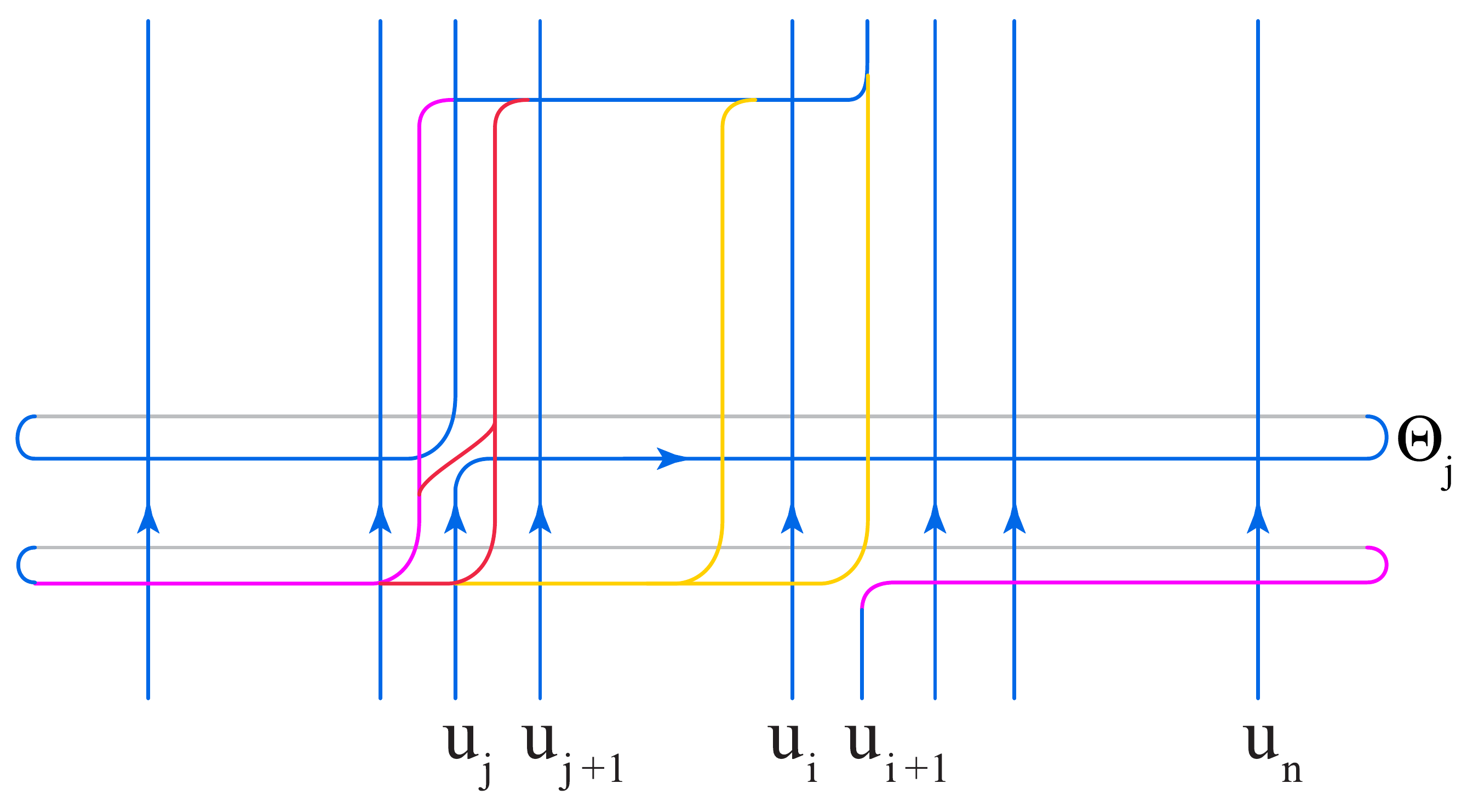}
\caption{}
\label{pic-deltaijmut2}
\end{figure}
Then we use the formula
\bea
 \Theta_i \Theta_j\RR_{\ov{j\, i}}=\RR_{\ov{j\, i}} \Theta_i \Theta_j.\nn
\eea
Hence we have the following expression
\bea
\Delta_i \Delta_j&=&\RR_{\ov{i\, i-1}}\ldots \RR_{\ov{i,j+1}}{\color{red}\RR_{\ov{i,j}}}{\color{blue}\RR_{\ov{i,j-1}}\ldots \RR_{\ov{i\, 1}}}{\color{magenta}\RR_{\ov{j\, j-1}}\ldots \RR_{\ov{j\, 1}}}{\color{green}\RR_{\ov{j\, i}}}\nn\\
&\times& \Theta_i \Theta_j\RR_{\ov{j\,n}}\ldots\RR_{\ov{j\, i+1}}
\RR_{\ov{j\,i-1}}\ldots\RR_{\ov{j\,j+1}}\RR_{\ov{i\,n}}\ldots\RR_{\ov{i\, i+1}}.\nn
\eea

Now, we perform the transposition of  blue and magenta elements with the help  of the  YBE involving the  red operator:
\bea
\Delta_i \Delta_j&=&{\color{magenta}\RR_{\ov{j\, j-1}}\ldots \RR_{\ov{j\, 1}}}\RR_{\ov{i\, i-1}}\ldots \RR_{\ov{i,j+1}}{\color{blue}\RR_{\ov{i,j-1}}\ldots \RR_{\ov{i\, 1}}}{\color{red}\RR_{\ov{i,j}}}{\color{green}\RR_{\ov{j\, i}}}\nn\\
&\times& \Theta_i \Theta_j\RR_{\ov{j\,n}}\ldots\RR_{\ov{j\, i+1}}
\RR_{\ov{j\,i-1}}\ldots\RR_{\ov{j\,j+1}}\RR_{\ov{i\,n}}\ldots\RR_{\ov{i\, i+1}}.\nn
\eea
By using the fact that
\bea
{\color{red}\RR_{\ov{i,j}}}{\color{green}\RR_{\ov{j\, i}}}={\color{green}\RR_{\ov{j\, i}}}{\color{red}\RR_{\ov{i,j}}}=1\nn
\eea
 we arrive to the formula
\bea
\Delta_i \Delta_j
&=&\RR_{\ov{j\, j-1}}\ldots \RR_{\ov{j\, 1}} \Theta_j\RR_{\ov{i\, i-1}}\ldots \RR_{\ov{i,j+1}}{\color{blue}\RR_{\ov{i,j-1}}\ldots \RR_{\ov{i\, 1}}}{\color{magenta}\RR_{\ov{j\,n}}\ldots\RR_{\ov{j\, i+1}}} \nn\\
&\times&
\RR_{\ov{j\,i-1}}\ldots\RR_{\ov{j\,j+1}} \Theta_i\RR_{\ov{i\,n}}\ldots\RR_{\ov{i\, i+1}}.\nn
\eea

Then transposing the commuting elements (blue ones and magenta ones) we finalize by the expression:
\bea
\Delta_i \Delta_j&=&
\RR_{\ov{j\, j-1}}\ldots \RR_{\ov{j\, 1}} \Theta_j{\color{magenta}\RR_{\ov{j\,n}}\ldots\RR_{\ov{j\, i+1}}} \RR_{\ov{i\, i-1}}\ldots \RR_{\ov{i,j+1}}\RR_{\ov{j\,i-1}}\ldots\RR_{\ov{j\,j+1}} \nn\\
&\times&
{\color{blue}\RR_{\ov{i,j-1}}\ldots \RR_{\ov{i\, 1}}} \Theta_i\RR_{\ov{i\,n}}\ldots\RR_{\ov{i\, i+1}}.\nn
\eea
Now, we insert the product 
\bea
{\color{green}\RR_{\ov{j\, i}}}{\color{red}\RR_{\ov{i,j}}}\nn
\eea
and by using the YBE several times in the same fashion as below we arrive to the following expression
\bea
\Delta_i \Delta_j&=&
\RR_{\ov{j\, j-1}}\ldots \RR_{\ov{j\, 1}} \Theta_j\RR_{\ov{j\,n}}\ldots\RR_{\ov{j\, i+1}}{\color{green}\RR_{\ov{j\, i}}} {\color{red}\RR_{\ov{i,j}}}{\color{blue}\RR_{\ov{i\, i-1}}\ldots \RR_{\ov{i,j+1}}}{\color{magenta}\RR_{\ov{j\,i-1}}\ldots\RR_{\ov{j\,j+1}}}\nn\\
&\times&
\RR_{\ov{i,j-1}}\ldots \RR_{\ov{i\, 1}} \Theta_i\RR_{\ov{i\,n}}\ldots\RR_{\ov{i\, i+1}}\nn\\
&=&
\RR_{\ov{j\, j-1}}\ldots \RR_{\ov{j\, 1}} \Theta_j\RR_{\ov{j\,n}}\ldots\RR_{\ov{j\, i+1}}{\color{green}\RR_{\ov{j\, i}}} {\color{magenta}\RR_{\ov{j\,i-1}}\ldots\RR_{\ov{j\,j+1}}}{\color{blue}\RR_{\ov{i\, i-1}}\ldots \RR_{\ov{i,j+1}}}{\color{red}\RR_{\ov{i,j}}}\nn\\
&\times&
\RR_{\ov{i,j-1}}\ldots \RR_{\ov{i\, 1}} \Theta_i\RR_{\ov{i\,n}}\ldots\RR_{\ov{i\, i+1}}=\Delta_j\Delta_i.\nn
\eea
 This completes the proof.

\end{document}